\documentclass{amsart}

\newtheorem{theorem}{Theorem}
\newcommand{\bt}{\begin{theorem}}
\newcommand{\et}{\end{theorem}}

\newtheorem*{theoremNN}{Theorem}
\newcommand{\btNN}{\begin{theoremNN}}
\newcommand{\etNN}{\end{theoremNN}}

\newtheorem{lemma}{Lemma}
\newcommand{\bl}{\begin{lemma}}
\newcommand{\el}{\end{lemma}}
\newtheorem{corollary}{Corollary}
\newcommand{\bc}{\begin{corollary}}
\newcommand{\ec}{\end{corollary}}
\newcommand{\beq}{\begin{equation}}
\newcommand{\eeq}{\end{equation}}
\newcommand{\benum}{\begin{enumerate}}
\newcommand{\eenum}{\end{enumerate}}

\newcommand{\R}{\ensuremath{\mathbf R}}

\newcommand{\bsmallmat}{\left(\begin{smallmatrix}}
\newcommand{\esmallmat}{\end{smallmatrix}\right)}

\newcommand{\bmat}{\left(\begin{matrix}}
\newcommand{\emat}{\end{matrix}\right)}
\DeclareMathOperator{\qand}{\quad\text{and}\quad}
\DeclareMathOperator{\qqand}{\qquad\text{and}\qquad}

\DeclareMathOperator{\sgn}{\text{sign}}

\newtheorem{problem}{Problem}
\newcommand{\bprob}{\begin{problem}}
\newcommand{\eprob}{\end{problem}}

\usepackage{amssymb,latexsym}

\title{Gauss's proof of Descartes's rule of signs}

\author{Melvyn B. Nathanson}

\address{Department of Mathematics\\Lehman College (CUNY)\\Bronx, NY 10468}

\email{melvyn.nathanson@lehman.cuny.edu}



\date{\today}

\begin{document}

\begin{abstract} 
For nonzero polynomials $f(x)$ with real coefficients, 
Descartes's rule of signs gives an upper bound for the number of positive zeros.  
A variation of an argument  of Gauss is used to prove Descartes's theorem. 
 \end{abstract}

\maketitle

\section{Variation of sign}

The \emph{sign}\index{sign} of a real number $\lambda$ is the function 
\[
\sgn(\lambda) = \begin{cases}
1 & \text{ if $\lambda > 0$} \\ 
0 & \text{ if $\lambda = 0$} \\ 
-1 & \text{ if $\lambda < 0$.} 
\end{cases}
\]
Let 
\[
A = \left( a_0,a_1,a_2,\ldots, a_n \right)
\]
be a sequence of real numbers.  
The associated \emph{sign sequence} is 
\[
\sgn(A) = \left( \sgn(a_0), \sgn(a_1), \sgn(a_2), \ldots, \sgn(a_n) \right). 
\]
Let $k+1$ be the number of nonzero terms in this sequence.  
There is a strictly increasing sequence $0 \leq j_0 < j_1 < \cdots < j_k \leq n$ such that  
$ \left( a_{j_0}, a_{j_1}, \ldots, a_{j_k} \right)$ is the subsequence of nonzero terms in $A$. 
We have $\sgn\left( a_{j_i}\right) = \pm 1$ for all $i \in \{0,1,\ldots, k\}$. 
The number of  \emph{sign variations}\index{sign variations} 
or \emph{sign changes}\index{sign changes} in the sequence $A$, denoted $V(A)$, 
 is the number of pairs $(a_{j_{i-1}}, a_{j_i})$ 
such that $\sgn(a_{j_{i-1}}) = -\sgn( a_{j_i})$.  
The number of sign changes in the sequence $(a_0,a_1,a_2,\ldots, a_n)$  
equals the number of sign changes in the associated sign sequence $\sgn(A)$. 

For example, if $A = (1,-1,0,1,-1)$, $A_1 = (1,-1,0) $, and  $A_2 =(0,1,-1)$, then 
\[
V(A) = 3, \quad V(A_1) = 1, \qqand V(A_2) = 1.
\]
Note that 
\[
V(A) \neq V(A_1) +  V(A_2). 
\]
Consider $n = 4$ and $\ell = 2$ in the following lemma.

\bl                                      \label{Descartes:lemma:additionFormula} 
Let  $A = \left( a_0,a_1,a_2,\ldots, a_n \right)$ 
be a sequence of real numbers.  
For $\ell \in \{1,2,\ldots, n-1\}$, let 
\[
A_1 = \left( a_0,a_1,\ldots, a_{\ell} \right) 
\]
and 
\[
A_2 = \left( a_{\ell},  a_{\ell+1}, \ldots, a_n \right). 
\]
If $a_{\ell} \neq 0$, then 
\[   
V(A) = V(A_1) + V(A_2). 
\]
If  $a_{i} = 0$ for all $i \in \{0,1,\dots, \ell\}$ 
or if $a_j = 0$ for all $j \in \{\ell,\ell+1,\ldots, n\}$,  
then 
\[   
V(A) = V(A_1) + V(A_2). 
\]
If $a_{\ell} = 0$ and if $a_i \neq 0$ for some $i < \ell$ and $a_j \neq 0$ for some $j > \ell$  
and if $\ell_1$ is the largest integer such that $\ell_1 < \ell$ and $a_{\ell_1} \neq 0$, 
and $\ell_2$ is the smallest integer such that $\ell_2 >\ell$ and $a_{\ell_2} \neq 0$, 
then
\[   
V(A) = 
\begin{cases}
V(A_1) + V(A_2) 
& \text{if  $\sgn(a_{\ell_1}) = \sgn(a_{\ell_2})$}\\
V(A_1) + V(A_2) + 1
& \text{if $\sgn(a_{\ell_1}) \neq \sgn(a_{\ell_2})$.}
\end{cases} 
\]
\el

\begin{proof}
This is a straightforward verification.  
\end{proof}

\bl                      \label{Descartes:lemma:Descartes-parity}
Consider a sequence 
$\varepsilon = ( \varepsilon_0,\varepsilon_1,\ldots, \varepsilon_k)$ 
such that  $\varepsilon_i \in \{1,0,-1\}$ for all $i \in \{0,1,\ldots, n\}$ and 
$\varepsilon_0\varepsilon_k \neq 0$. 
\benum
\item[(i)]
If $\varepsilon_0 = \varepsilon_k$, then $V(\varepsilon)$ is even.  
\item[(ii)]
If $\varepsilon_0 = -\varepsilon_k$, then $V(\varepsilon)$ is odd. 
\item[(iii)]
If $\varepsilon_{i-1}\varepsilon_{i+1} = -1$ for some $i \in \{1,\ldots, k-1\}$,  then 
$V(\varepsilon_{i-1}, \varepsilon_{i}, \varepsilon_{i+1}) = 1$ for all $\varepsilon_i \in \{0,1,-1\}$. 
\eenum
\el

\begin{proof}
Statements~(i) and~(ii)  follow from the simple but fundamental observation 
that changing signs an even number of times does not change the sign 
but changing signs an odd number of times does change the sign.

To prove statement~(iii), we observe that if $\varepsilon_{i-1} = -\varepsilon_{i+1} = 1$, 
then $(\varepsilon_{i-1}, \varepsilon_{i}, \varepsilon_{i+1}) = (1,1,-1)$ or $(1,0,-1)$ or 
$(1,-1,-1)$ and each of these sequences has one sign change.  
Similarly, if  $\varepsilon_{i-1} = -\varepsilon_{i+1} = -1$, 
then $(\varepsilon_{i-1}, \varepsilon_{i}, \varepsilon_{i+1}) = (-1,1,1)$ or $(-1,0,1)$ or 
$(-1,-1,1)$ and each of these sequences has one sign change.  
This completes the proof.  
\end{proof}

A $C^n$ function is a  function with $n$  continuous derivatives.  
Let $f(x)$ be a $C^n$ function.  For $\lambda$ in the domain of $f(x)$,  
we define the \emph{derivative sequence}\index{derivative sequence} 
\[
D_f^{(n)}(\lambda) = \left( f(\lambda), f'(\lambda), f''(\lambda), \ldots, f^{(n)}(\lambda) \right).
\]
The variation of sign  function  
\[
V^{(n)}_f(\lambda) = V\left(D_f^{(n)}(\lambda)  \right) 
= V\left(  f(\lambda), f'(\lambda), f''(\lambda), \ldots, f^{(n)}(\lambda)  \right)
\]
counts the number of sign changes in the derivative sequence.

In this paper we consider only polynomials with real coefficients.  
Let $f(x)$ be a polynomial of degree $n$.  For all $\lambda \in \R$, 
there is the Taylor expansion 
\[
f(x) = \sum_{j=0}^n a_j (x-\lambda)^j
\]
with Taylor coefficients  
\[
a_j= \frac{f^{(j)}(\lambda)}{j!}.  
\]
It follows that   
\[
\sgn(a_j) = \sgn\left( f^{(j)}(\lambda)\right) 
\]
for all $j \in \{ 0,1,\ldots, n\}$ and so 
\[
V_f^{(n)}(\lambda) = V\left( a_0, a_1, \ldots, a_n \right)
\]
 also counts the number of sign changes in the sequence of 
Taylor coefficients of $f(x)$.

For example, the polynomial 
\[
f(x) = 10 + 8x  - 3x^2   - 5x^3  + 2x^5  +7x^8 + x^9
\]
has coefficient sequence  
\[
A = (10, 8, -3, -5, 0, 2, 0, 0, 7, 1). 
\]  
The corresponding sign sequence is 
\[
\sgn(A) = (1, 1, -1, -1, 0, 1, 0, 0, 1, 1).
\]
There are two sign changes in these sequences and so $V^{(n)}_f(0) = 2$.

Let $f(x)$ be a $C^n$ function and let 
\[
g(x) = f(x+a).
\]
For all $j \in \{0,1,\ldots, n\}$ we have 
$g^{(j)}(x) = f^{(j)}(x+a)$ and so 
\[
g^{(j)}(0) = f^{(j)}(a). 
\]
Therefore, 
\begin{align}            \label{Descartes-Vg}
V_g^{(n)}(0) & = V\left(g(0), g'(0), g''(0), \ldots, g^{(n)}(0) \right)        \nonumber  \\
& = V\left(f(a), f'(a), f''(a), \ldots, f^{(n)}(a) \right) \\
& = V_f^{(n)}(a).      \nonumber
\end{align}

\section{Counting zeros} 
Let $f(x)$ be a nonzero polynomial with real coefficients. 
Let $\mu_f(\lambda)$ be the multiplicity of  the real or complex number $\lambda$  
as a root of the polynomial $f(x)$.  We define  the \emph{zero counting function} 
\[
Z_f(a,b) =  \sum_{a < \lambda < b}  \mu_f(\lambda). 
\]
This is  the number  (counting multiplicity) of real numbers $\lambda$ such that 
\[
f(\lambda) = 0 \qqand a < \lambda < b.
\]
In particular, $Z_f(0,\infty)$ counts the number of positive roots of $f(x)$ 
and $Z_f(-\infty,0)$ counts the number of negative roots of $f(x)$.  

We define similarly 
\[
Z_f(a,b] =  \sum_{a < \lambda \leq b}  \mu_f(\lambda). 
\qand
Z_f[a,b) =  \sum_{a \leq \lambda < b}  \mu_f(\lambda). 
\]

Let $f(x)$ be a nonzero polynomial and let $g(x) = f(-x)$.  
For $\lambda > 0$ we have $f(-\lambda) = 0$ if and only if 
$g(\lambda) = 0$, and so 
\[
Z_f(-\infty, 0) = Z_g(0,\infty). 
\]
For $a,b \in \R$ with $a < b$, we have  
\beq            \label{Descartes-5}
Z_f(a,b] = Z_f(a,\infty) - Z_f(b,\infty).
\eeq

Let $g(x) = f(x+a)$.  
Equivalently, $f(x) = g(x-a)$.  
We have $f(\lambda) = 0$ if and only if $g(\lambda - a) = 0$.  
Because $\lambda > a$ if and only if $\lambda - a > 0$, 
it follows that 
\[
Z_f(a,\infty) = Z_g(0,\infty).
\]
Similarly, if $h(x) = f(x+b)$, then 
\[
Z_f(b,\infty) = Z_h(0,\infty).
\]
Relation~\eqref{Descartes-5} implies  that if $f(b) \neq 0$, then 
\beq            \label{Descartes-6}
Z_f(a,b]  = Z_g(0,\infty) - Z_h(0,\infty).  
\eeq 
We may use  Descartes's rule of signs (Theorem~\ref{Descartes:theorem:Descartes})  
to compute zero counting 
functions of the form $Z_f(a,b]$.

If $\mu = \mu_f(0)$, 
then $f(x) = x^{\mu}\ell(x)$, where $\ell(x)$ is a polynomial of degree $n - \mu$ 
such that $\ell(0) \neq 0$.  
For $\lambda \neq 0$, we have $f(\lambda) = 0$ if and only if $\ell(\lambda) = 0$,
and so  the nonzero roots of $f(x)$ are the nonzero roots of $\ell(x)$.
Moreover, $\mu_f(\lambda) = \mu_\ell(\lambda)$ for all $\lambda \neq 0$. 
Thus, to count the number of nonzero roots of  polynomials, it suffices to
consider only polynomials $f(x)$ with $f(0) \neq 0$, that is, 
polynomials with nonzero constant terms.

\bt[Descartes]                   \label{Descartes:theorem:Descartes}
Let $f(x)$ be a  polynomial with $f(0) \neq 0$.   There exists 
a nonnegative even integer $\nu$ such that  
\[
Z_f(0,\infty) = V^{(n)}_f(0) - \nu. 
\] 
\et

\bc
If $V^{(n)}_f(0) = 1$, then $f$ has exactly one positive root $\lambda$ and $\lambda$ is a simple root.  
\ec

\begin{proof}
Because $\nu$ is nonnegative and even, if $\nu \neq 0$, then $\nu \geq 2$ and 
\[
0 \leq Z_f(0,\infty) = V^{(n)}_f(0) - \nu   = 1 - \nu \leq -1
\]
which is absurd.  Therefore, $\nu = 0$ and $ Z_f(0,\infty) = V^{(n)}_f(0) - \nu = 1$. 
This completes the proof.  
\end{proof}

Let $a,b \in \R$ with $a < b$.  
Let $f(x)$ be a  polynomial of degree $n$ with $f(b) \neq 0$.  
Consider the polynomials $g(x) = f(x+a)$ and $h(x) = f(x+b)$. 
 Descartes's rule of signs and relation~\eqref{Descartes-6} give 
 nonnegative even  integers $\nu_g$ and $\nu_h$ such that 
\[
Z_f(a,\infty) = Z_g(0,\infty) = V^{(n)}_g(0) - \nu_g = V_f^{(n)}(a) - \nu_g  
\]
\[
Z_f(b,\infty) = Z_h(0,\infty) = V^{(n)}_h(0) - \nu_h = V_f^{(n)}(b) - \nu_h  
\]
From equation~\eqref{Descartes-6} we obtain   
\[
Z_f(a,b] =  V_f^{(n)}(a) - V_f^{(n)}(b) - \nu
\]
where 
\[
\nu = \nu_g - \nu_h.
\]

This argument does not prove the nonnegativity of the integer $\nu$.
The nonnegativity of $\nu$ would immediately imply the following result, 
which is the Budan-Fourier theorem. 

\bt [Budan-Fourier]    \label{realrooted:theorem:Budan-Fourier}
Let $f(x)$ be a polynomial of degree $n$ and let $a,b \in \R$ with $a<b$.    
There exists a nonnegative even integer $\nu$ such that  
\[
Z_f(a,b) =  V^{(n)}_f(a) - V^{(n)}_f(b) - \nu.
\]
\et

 The Budan-Fourier theorem implies Descartes's rule of signs.

For the history and alternate proofs of Descartes's theorem, 
see~\cite{basu-poll-roy06, bens13, conk43, dick22, jaco74, hurw12}.  
A recent proof of the Budan-Fourier theorem is Nathanson~\cite{nath22z}.

\section{Gauss's proof of Descartes's theorem} 

The following result is the essential part of Gauss's beautiful  
proof~\cite{gaus28} of Descartes's theorem.

\bl                          \label{realrooted:lemma:Descartes-Gauss}
Let $h(x)$ be a  nonzero polynomial with $h(0) \neq 0$ 
and let $\lambda > 0$.  If  
\[
\ell(x) = (x-\lambda)h(x)
\]
then 
\[
V^{(n)}_{\ell}(0) =  V^{(n)}_h(0) +1 +  \nu
\] 
for some nonnegative even integer $\nu$.  
\el

\begin{proof}
Let 
\[
h(x)  =  \sum_{j=0}^{m} c_j x^j 
\]
be a polynomial of degree $m$ with $h(0) = c_0 \neq 0$ and let 
\[
V^{(m)}_h(0) = V(c_0, c_1,\ldots, c_m) = v 
\] 
be the number of sign changes in the coefficient sequence of $h(x)$.  
Note that 
\[
 v \in \{0,1,\ldots,m\}. 
 \]
We shall construct an integer sequence $\left( j_i\right)_{i=0}^v$  
with 
\beq            \label{Descartes-0}
0 = j_0 < j_1 < j_2 < \cdots < j_v \leq m 
\eeq
such that 
\beq            \label{Descartes-1}
V(c_{j_{i-1}}, \ldots, c_{j_i - 1}) = 0
\qqand 
V(c_{j_{i-1}}, \ldots, c_{j_i}) = 1
\eeq
for all $i = 1,2,\ldots, v$, and 
\beq            \label{Descartes-2}
V( c_{j_v},\ldots, c_m) = 0. 
\eeq

The construction is by induction. 
Let $j_0 = 0$.  Let $i \in \{1,\ldots, v\}$. 
If $j_{i-1} \in \{0,1,\ldots, m\}$ satisfies  
\[
V(c_0,\ldots, c_{j_{i-1}}) = i-1 \leq v-1 < V^{(m)}_h(0)
\]
then there is a smallest integer $j_i > j_{i-1}$ such that 
\[
\sgn(c_{j_i}) = -\sgn(c_{j_{i-1}}) 
\] 
and so 
\[
V(c_0,\ldots, c_{j_{i-1}}, c_{j_i} ) = i. 
\]
The sequence $(j_i)_{i=0}^v$ satisfies conditions~\eqref{Descartes-0} --~\eqref{Descartes-3}.
Moreover, if $j_v < m$, then 
\beq            \label{Descartes-4}
\sgn(c_{j_v})  =\sgn(c_m). 
\eeq

Consider the polynomial 
\begin{align*}
\ell(x) & = (x-\lambda)h(x) 
 =  \sum_{j=0}^{m} c_j x^{j+1}  -  \sum_{j=0}^{m} \lambda c_j x^j \\
& =  - \lambda c_0 + \sum_{j=1}^{m} \left( c_{j-1} - \lambda c_j \right) x^j   + c_m x^{m+1}   \\
& =  \sum_{j=0}^{m+1} b_j x^j  
\end{align*}
with 
\[
b_j = \begin{cases} 
 -\lambda c_0 & \text{for $j = 0$} \\
c_{j-1} -\lambda c_j & \text{for $j = 1, \ldots, m$} \\
c_m & \text{for $j = m+1$.} 
\end{cases} 
\]
For all $i \in \{1,\ldots, v\}$ we have 
\[
b_{j_i} = c_{j_i-1} -\lambda c_{j_i}. 
\]
It follows from the construction of the sequence $(j_i)_{i=0}^v$ 
that either 
\[
c_{j_i-1} = 0  
\]
or 
\[
 \sgn(c_{j_i-1})  =   \sgn(c_{j_{i-1}})  = -  \sgn(c_{j_i}).  
 \]
Because $\lambda > 0$, we have 
\[
\sgn(b_{j_i}) = \sgn(c_{j_i-1} - \lambda c_{j_i}) = -\sgn(c_{j_i})  
\]
and so  
\[
\sgn(b_{j_i}) = -\sgn(c_{j_i})  = \sgn(c_{j_{i-1}})  = - \sgn(b_{j_{i-1}}).   
\]
We have  $j_v \leq m$.  It follows from~\eqref{Descartes-4} that  
\[
\sgn(b_{j_v}) = -\sgn(c_{j_v})  = -\sgn(c_{m}) = -\sgn(b_{m+1}).  
\]\
By Lemma~\ref{Descartes:lemma:Descartes-parity}, 
there exist nonnegative even integers $\nu_0, \nu_1,\ldots, \nu_v$ 
such that 
\[
V(b_{j_{i-1}}, \ldots, b_{j_i} ) = 1 + \nu_{i-1}
\]
 for all $i \in \{1,\ldots, v\}$ and  
\[
V(b_{j_v}, \ldots, b_{m+1} ) = 1 + \nu_{v}. 
\]   
Because $b_{j_i} \neq 0$ for all $i \in \{1,\ldots, v\}$, 
the addition formula in Lemma~\ref{Descartes:lemma:additionFormula} gives  
\begin{align*}
V^{(m+1)}_{\ell}(0)
&  = V(b_0,\ldots, b_{m+1} ) = \sum_{i=1}^v V( b_{j_{i-1}}, \ldots, b_{j_i} ) 
+ V(b_{j_v}, \ldots, b_{m+1} ) \\
&  =  v + 1 +\sum_{i=0}^{v} \nu_i \\
&  =  V^{(m)}_h(0)  + 1 + \nu  
\end{align*} 
where $\nu = \sum_{i=0}^{v} \nu_i$ is a nonnegative even integer. 
This completes the proof. 
\end{proof}

\bl          \label{realrooted:lemma:Descartes-5}
Let $h(x)$ be a nonzero polynomial with $h(0) \neq 0$.  
Let $(\lambda_i)_{i=1}^k$ be a sequence of not necessarily distinct positive numbers. 
If 
\[
f(x) =   h(x) \prod_{i=1}^k (x-\lambda_i) 
\]
then
\[
V^{(n)}_f(0)= k + V^{(n)}_h(0) + \nu
\]
for some nonnegative even integer $\nu$.  
\el

\begin{proof}
This follows from Lemma~\ref{realrooted:lemma:Descartes-Gauss}
by induction on $k$.
\end{proof}

We now prove Descartes's rule of signs (Theorem~\ref{Descartes:theorem:Descartes}). 

\begin{proof}
Let $f(x)$ be a monic polynomial of degree $n$ with $f(0) \neq 0$ and let  
\[
Z_f(0,\infty) = k.
\]

Let $\Re(z)$ denote the real part of the complex number $z$. 

Let $(\lambda_i)_{i=1}^n$ be the sequence of $n$ not necessarily distinct 
real and complex  roots of $f(x)$.  Note that $\lambda_i \neq 0$ for all $i \in \{1,\ldots, n\}$. 
We  order the roots so that, 
for nonnegative integers $k,\ell, m$ with 
\[
k+\ell + 2m = n
\]
we have  
\benum
\item[(i)]
$(\lambda_i)_{i=1}^k$ is the sequence of positive roots of $f(x)$,
\item[(ii)]
$(\lambda_i)_{i=k+1}^{k+\ell}$ is the sequence of negative roots of $f(x)$, 
\item[(iii)]
$(\lambda_i)_{i= k+\ell+1}^{k+\ell + 2m}$ is the sequence of nonreal complex  
roots of $f(x)$ with  
\[
\overline{\lambda_{k+\ell+r}} = \lambda_{k+\ell+r+m}. \qquad\text{for $r\in \{1,\ldots, m\}$.}
 \]
\eenum
For every nonreal complex number $\lambda_i$, we have 
\[
(x-\lambda_i)(x-\overline{\lambda_i}) = x^2 -2\Re(\lambda_i)x + |\lambda_i|^2. 
\]
The constant term of this quadratic polynomial  is positive 
and so the constant term of the polynomial 
\begin{align*} 
\prod_{i=k+\ell+1}^{k+\ell+2m} (x-\lambda_i) 
& =  \prod_{i=k+\ell+1}^{k+\ell+m} (x-\lambda_i)(x-\overline{\lambda_i}) \\ 
& =  \prod_{i=k+\ell+1}^{k+\ell+m} \left( x^2 -2\Re(\lambda_i)x + |\lambda_i |^2 \right) 
\end{align*} 
is positive:
\[
 \prod_{i=k+\ell+1}^{k+\ell+m}  |\lambda_i |^2 > 0
\]
Because $\lambda_i < 0$ for all $i \in \{k+1,\ldots, k+\ell\}$, 
the constant term of the polynomial 
\[
\prod_{i=k+1}^{k+\ell} (x-\lambda_i) 
\]
is also positive: 
\[
 \prod_{i=k+1}^{k+\ell}  |\lambda_i|  > 0.
\]
Therefore, the constant term of the monic polynomial 
\[
h(x)  = \prod_{i=k+1}^n (x-\lambda_i)
\]
is positive: 
\[
h(0) = \prod_{i=k+1}^{k+\ell}  |\lambda_i| \prod_{i=k+\ell+1}^{k+\ell+m}  |\lambda_i |^2 > 0.  
\]
The leading coefficient and the constant term of $h(x)$ are positive.   
By Lemma~\ref{Descartes:lemma:Descartes-parity}, 
the number of sign changes  in the coefficient sequence of $h(x)$ 
is the nonnegative even integer $V^{(n - k)}_h(0) = \nu_h$. 

We have 
\[
f(x) =   h(x)  \prod_{i=1}^k (x - \lambda_i).  
\]
Because $(\lambda_i)_{i=1}^k$ is a sequence of positive real numbers,  Lemma~\ref{realrooted:lemma:Descartes-5} implies that 
there is a nonnegative  even integer $\nu_f$ such that 
\[
V^{(n)}_f(0)=  V^{(n)}_h(0) + k + nu_f = k + \nu_h + \nu_f 
= k +  \nu
\]
for $\nu = \nu_h + \nu_f$, and so 
\[
Z_f(0,\infty) = k = V^{(n)}_f(0) - \nu. 
\]
This completes the proof. 
\end{proof}

\end{document}